\documentclass[12pt,reqno]{article}
\def\hybrid{\topmargin 0pt      \oddsidemargin 0pt
        \headheight 0pt \headsep 0pt
        \textwidth 170true mm       % US paper
        \textheight 231true mm         % US paper
        \marginparwidth 0.0in
        \parskip 0pt plus 1pt   \jot = 1.5ex}

\usepackage{amssymb}
\usepackage{amsmath}
\usepackage{amsthm}
\hybrid
\usepackage{amssymb}
\usepackage{amsmath}
\usepackage{amsthm}

\newcommand{\kf}{{\cal K}_{\alpha}^{(1)}}
\newcommand{\ks}{{\cal K}_{\alpha}^{(2)}}

\renewcommand{\Re}{\mathrm{Re}\;}
\renewcommand{\Im}{\mathrm{Im}\;}
\newcommand{\K}{K^\alpha}
\newcommand{\me}{m_{\alpha,\varepsilon}(|\xi|)}
\newcommand{\F}{ {\;\;_{2}F_{1}}}
\newcommand{\R}{{\Bbb R}^n}
\newcommand{\Sa}{S^\alpha}
\newcommand{\Na}{N^\alpha}
\newcommand{\N}{N^\alpha}
\newcommand{\normf}{\|f\|_{p}}
\newcommand{\C}{\Bbb C}

\newcommand{\lp}{L_{p}}
\newcommand{\lr}{L_{q}}
\newcommand{\ma}{m_\alpha(|\xi|)}

\newcommand{\intR}{\int\limits_{{\Bbb R}^n}}
\newcommand{\Sal}{S^\alpha_\ell}

\newtheorem{remark}{Remark}[section]
\newtheorem{lemma}{Lemma}[section]
\newtheorem{theorem}{Theorem}[section]

\title{$(\lp, \lr)$ estimates of potentials with oscillating
kernel\\[18pt]}
\author{E. Ournycheva }
\date{}
\begin{document}

\thispagestyle{empty} \maketitle

\begin{center}
Department of Mathematics, Bar-Ilan University, 52900 Ramat Gan,
Israel\\
 e-mail: ournyce@macs.biu.ac.il\\[18pt]
\end{center}

\begin{abstract}
\noindent
$(L_p,  L_q)$ estimates are obtained for oscillatory  potentials
\begin{equation}
(K^\alpha
f)(x)=\int\limits_{R^n}\frac{\exp(i|y|)}{|y|^{n-\alpha}}f(x-y)dy,
\quad 0<\alpha<n, \quad   n\geq 2, \nonumber \end{equation} whose
symbol has a singularity on the unit sphere. These potentials are
natural modifications of the celebrated  Bochner-Riesz operator
and  Helmholtz potential arising  in  Fourier analysis and PDE.
For some values of $\alpha$, determination of the corresponding
pairs $(L_p, L_q)$ represents an open problem. The range of
$\alpha$ for which the problem is open is just the same as for the
Bochner-Riesz means.
\end{abstract}

\setcounter{equation}{0} \setcounter{theorem}{0}
\setcounter{lemma}{0}
\section{Introduction}

We consider the  potential operator with oscillating kernel

\begin{equation}
(\K f)(x)=\intR\frac{\exp(i|y|)}{|y|^{n-\alpha}}f(x-y)dy, \quad 0<\alpha<n, \quad n\geq 2.
\label{0.1}
\end{equation}
Up to now,  the  $(\lp, \lr)$ estimates for oscillatory potentials
on $\R$ have been investigated  for  some special cases only. We
mention the oscillation generated by the Bessel function
 (Bochner-Riesz means) and  that  generated by the Hankel
function  (Helmholtz potential). Such oscillatory integrals  arise from the classical
 problems of Fourier analysis concerning  the summability of
series and inversion of the Fourier transform  and from the
Helmholtz equation with the Dirichlet boundary condition.
Respectively, more general potentials (\ref{0.1}) are natural
modification of those operators.

The multiplier problem for the Bochner-Riesz means has a long
history. We refer to \cite{Carleson,Cor, Davis, Fefferman1,
Fefferman2, H, Stein1} for the background information.
 We also point out the papers \cite{Bak, BMO, Borjeson, S} and \cite{Nogin, Rubin}
 where one can find $(\lp, \lr)$
estimates for the Bochner-Riesz operators of negative order and the Helmholtz potentials,
respectively.

For different $\alpha$, we establish the boundedness of $\K$
taking $\lp$ into $\lr$. The sets of pairs
$(\frac{1}{p},\frac{1}{q})$  for which the operator $K^\alpha$ is
bounded from $L_p$ into $L_q$ are convex and of the special form
on the $(\frac{1}{p},\frac{1}{q})$-plane (see Figure). The $\cal
L$-characteristic of the operator $\K$ is constructed for either
$n=2$, or $n>2$, provided $0<\alpha<\frac{n(n-1)}{2(n+1)}$ or
$\frac{n}{2}\leq\alpha<n$. In other cases
 we  have  gaps between necessary and sufficient  conditions
of boundedness similar to those   for the Bochner-Riesz operator.
Note, that the principal difficulties in the oscillatory
potentials theory are related to  controlling the oscillation of
kernels.
%Certain $(\lp, \lr)$ estimates for the operator (\ref{0.1}) have been
%obtained in \cite{}. We give some clarifications and improve the results of this paper.
The techniques used here to get through those difficulties goes
back to Stein \cite{Stein1}, Fefferman  \cite{Fefferman1}, Bak
\cite{Bak, BMO} and B\"orjeson \cite{Borjeson}.

The paper is organized as follows. In Section~\ref{sect1} we
formulate our main result (Theorem~\ref{t1}). Section~\ref{sect2}
contains  necessary preliminaries that will be used throughout
this paper.
 Section~\ref{sect3} deals with  some auxiliary statements.
First, the symbol $\ma$ of the operator $\K$ is estimated. Since
the symbol has a singularity on the unit sphere,  the operator
$\K$ shares properties of
 both the Bochner-Riesz means and
Riesz potential. More precisely, if we split the integral in
(\ref{0.1}) into two over $|y|\geq 1$ and $|y|< 1$ and
\begin{equation}\label{0.2}
\K f= \Sa f+\N f,
\end{equation}
where
\begin{equation}\label{4.6.1}
 (S^\alpha f)(x)=\int\limits_{|y|\geq
 1}\frac{\exp(i|y|)}{|y|^{n-\alpha}}f(x-y)dy,
\end{equation}

\begin{equation}\label{4.6.2}
(\N f)(x)=\int\limits_{|y|<
1}\frac{\exp(i|y|)}{|y|^{n-\alpha}}f(x-y)dy,
\end{equation}
then $\Sa$ will be responsible for the Bochner-Riesz means
properties while $\Na$  for properties of the Riesz potential.
Furthermore, we give estimates of an auxiliary operator, arising
from Stein's decomposition of the operator $\Sa$ (see
(\ref{4.7})). We next establish $\lp-$ boundedness of  $\K$.
Clearly, it can be reduced to the $\lp-$ boundedness property of
the operator $\Sa$. So the situation here is just the same as in
the Bochner-Riesz case, and the problem of precise $p$ for which
$\K$ is bounded on $\lp$ is correspondingly for $n>2$,
$\frac{n(n-1)}{2(n+1)}\leq\alpha<\frac{n-1}{2}$ (the problem of
Stein). In Subsection~\ref{s4.4} we give $(\lp, \lr)$ estimates of
 $\K$ along a segment
 through $(\frac{1}{p}, \frac{1}{p'})$, where $p'=p/(p-1)$, and perpendicular to the line of duality
 $\frac{1}{p}+ \frac{1}{q}=1$ provided
$\frac{n-1}{2}<\alpha <n$. We consider the cases $n=2$ and $n>2$
separately. In the first one the  result obtained here is sharp.
%The proof is based on  Bak's method (\cite{Bak}) using the Lorentz
%space estimate as a key point.
 When $n>2$, the result is sharp for  $n/2\leq\alpha<n$ only,  as in the case of
Bochner-Riesz operator with negative index.
% we apply complex interpolation with respect to $\alpha$.,
 The proof of Theorem~\ref{t1} is
given in Section~\ref{sect4}. For the sake of convenience, we give
one technical estimate in Appendix.

\setcounter{equation}{0}
\setcounter{theorem}{0}
\setcounter{lemma}{0}
\section{The main result}\label{sect1}

To formulate our main result we label some points in $Q=[0;1]\times [0;1]$ (see Figure):\\
 $L(\frac{1}{2};\frac{1}{2})$, $E(1;0)$, $F(1;1)$,
$A(1;1-\frac{\alpha}{n})$, $A'(\frac{\alpha}{n}; 0)$,
$H(1-\frac{\alpha}{n};1-\frac{\alpha}{n})$, $H'(\frac{\alpha}{n};\frac{\alpha}{n})$,
$C(\frac{3}{2}-\frac{2\alpha}{n-1};\frac{3}{2}-\frac{2\alpha}{n-1})$,
$C'(\frac{2\alpha}{n-1}-\frac{1}{2};\frac{2\alpha}{n-1}-\frac{1}{2})$,
$D(\frac{\alpha+1}{n+1};\frac{n-\alpha}{n+1})$,
$B(1-\frac{(n-1)(n-\alpha)}{n(n+1)};1-\frac{\alpha}{n})$,
$B'(\frac{\alpha}{n};\frac{(n-1)(n-\alpha)}{n(n+1)})$,
$G(\frac{n+3}{2n+2};1-\frac{\alpha}{n})$,
$G'(\frac{\alpha}{n};\frac{n-1}{2n+2})$,
$P(\frac{4\alpha-n+3}{2n+2};\frac{1}{2})$,
$P'(\frac{1}{2};\frac{3n-4\alpha-1}{2n+2})$.\\
Let  $(ABC\dots P)$ be an open polygon in $Q$ with vertices
$A,B,C,\dots P$, and let  $[ABC\dots P]$ be its closure. Let also
$(AB)$ be an open interval in $Q$ with the ends $A$ and $B$,
$[AB]$, $[AB)$, $(AB]$, be the closed and half open intervals,
respectively. For linear operator $\cal A$ defined on $\lp$ spaces
we denote ${\cal L(A)}=\{(\frac{1}{p},\frac{1}{q})\in Q: \|{\cal
A}\|_{\lp \to \lr}<\infty\}$.\\
We are now in a position  to formulate our results.

\begin{theorem}\label{t1}
I. The following embeddings are valid: \\
 1) If $0<\alpha<\frac{n(n-1)}{2(n+1)}$, $n>2$ or  $0<\alpha<\frac{1}{2}$, $n=2$ then
  $(A'H'HA)\cup(A'A)\cup(H'H)\subset{\cal L}(\K)$;\\
 2) If  $\frac{n(n-1)}{2(n+1)}\leq\alpha<\frac{n-1}{2}$, $n>2$ then
 $(A'G'C'CGA)\cup(A'A)\cup(C'C)\subset{\cal L}(\K)$;\\
 3) If $\frac{1}{2}<\alpha<2$, $n=2$ or $\frac{n}{2}\leq\alpha<n$, $n>2$ then
$(A'B'BA)\cup(A'A)\cup (BB')\subset{\cal L}(\K)$;\\
 4) If $\frac{n-1}{2}<\alpha<\frac{n}{2}$, $n>2$ then
$(A'G'P'PGA)\cup(A'A)\cup (P'P)\subset{\cal L}(\K)$;\\
 5) If $\alpha=\frac{1}{2}$, $n=2$ then $(A'B'BA)\cup(A'A)\subset{\cal L}(\K)$;\\
 6) If $\alpha=\frac{n-1}{2}$, $n>2$ then $(A'G'LGA)\cup(A'A)\subset{\cal L}(\K)$.\\
II. The operator $\K$ is unbounded from $\lp$ into $\lr$ whenever \\
 1) $(\frac{1}{p}, \frac{1}{q})\in [HAF]\cup[H'A'O]$; \\
 2) $(\frac{1}{p}, \frac{1}{q})\in [A'AE]\setminus (AA')$; \\
 3) $(\frac{1}{p}, \frac{1}{q})\in (BB'HH')$ provided
$\frac{n-1}{2}<\alpha<n$.
\end{theorem}

We wish here to formalize a part of the reasoning presented in the
proof of the positive results. As it was mentioned above, the idea
of splitting $\K$ into (\ref{0.2}) leads us to controlling the
oscillation of the kernel of  $\Sa$. For this goal, we apply a
standard model to get control of  operators whose symbol has a
singularity on the unit sphere. The method is based on the
decomposition
\begin{equation}\label{new:4.2}
S^\alpha=\sum\limits_{\ell=0}^\infty S^\alpha_\ell
\end{equation}
into the sum of operators with kernels supported on dyadic annuli.
In fact, we use Stein's decomposition (\cite{Stein1}) setting
\begin{equation}\label{4.7}
(S^\alpha_\ell f)(x)= 2^{(\alpha-n)\ell}\int\limits_{|y|\geq
1}\exp(i|y|)\psi(y/2^\ell)f(x-y)dy,
\end{equation}
where $\psi(y)=|y|^{\alpha-n}[\eta(y)-\eta(2y)]$ is a smooth
function supported in $\frac{1}{2}<|y|<2$, $\eta(y)=1$ as $|y|\leq
1$ and $\eta(y)=0$ as $|y|\geq 2$.

\begin{figure}[htb]
\begin{center}
\begin{picture}(310,310)

\put(10,10){\vector(1,0){300}}
\put(10,10){\vector(0,1){300}}
 \put(310,0){$\frac{1}{p}$}
\put(0,310){$\frac{1}{q}$}
\put(250,10){\line(0,1){240}}
\put(10,250){\line(1,0){240}}
 \put(230,15){$E$}
 \put(240,252){$F$}
\put(65,0){$A'$}
 \put(65,80){$H'$}
\put(0,0){$O$}
\put(250,185){$A$}
 \put(170,185){$H$}
 \put(10,10){\line(1,1){240}}
\put(130,130){\line(2,1){97}}
 \put(130,130){\line(-1,-2){50}}
 \put(80,10){\line(0,1){70}}
\put(180,180){\line(1,0){70}} \put(65,58){$B'$} \put(192,185){$B$}
\put(126,95){$P'$} \put(152,123){$P$} \put(65,25){$G'$}
\put(225,185){$G$} \put(80,60){\line(1,1){119}} \put(132,121){$D$}
\put(81,30){\line(3,5){46}} \put(230,179){\line(-5,-3){77}}
\put(95,110){$C'$} \put(140,155){$C$} \put(120, 135){$L$}
\put(81,30){\line(1,3){26}} \put(230,179){\line(-3,-1){75}}
\put(80,10){\line(1,1){170}}
% \put(80,40){\line(2,3){39}}
%\put(220,180){\line(-3,-2){56}}
% \put(120,310){$Figure$}
\end{picture}
\end{center}
\center{Figure}
\end{figure}

Then the study of the operator $\Sa$ can be reduced to the study
of the oscillatory operator $G_\lambda$ (see (\ref{3.2.1})).
$\lp-$ boundedness property and $(\lp,\lr)$ estimates of the
latter ( Lemma \ref{l5.1} and \ref{l5.2}) have been investigated
in \cite{Stein1} and \cite{S}, respectively (see also \cite{Bak},
\cite{Ho}). Such argument is used to prove $\lp-$ boundedness of
the operator $\Sa$ and $(\lp, \lr)$ estimates on the open segment
$(BB')$ in the case $n=2$. It is to be noted that the method of
proof of Lemma~\ref{l4.4.1} which involves the estimates for
characteristic functions $f=\chi_E$ implying the Lorentz space
estimate is not new (see e.g. \cite {Bak}).

 To obtain $(\lp, \lr)$ estimates of
 $\K$  along a segment
 through $(\frac{1}{p}, \frac{1}{p'})$ and perpendicular to the line of
 duality when $n>2$, we apply  the modification of an interpolation theorem for analytic
families of operators (Theorem \ref{t3.2}). A special case of
Stein's interpolation theorem (\cite{Stein2}) established in
\cite{BMO} need a slight modification for our case. Such approach
developed in \cite{BMO} and adopted here allows to obtain
$(\lp,\lr)-$ boundedness of $\K$ not only for the point $D$, but
for indices $(\frac{1}{p}, \frac{1}{q})$ off the line of duality.
Altogether then, it is relatively simple matter to construct the
$\cal L-$ characteristic of $\K$ in two dimension as well as in
higher dimensions provided   $0<\alpha<\frac{n(n-1)}{2(n+1)}$ or
$\frac{n}{2}\leq\alpha<n$. In other cases we establish additional
estimates for the operators $\Sal$. The arguments that are
required for proofs are based on properties of the symbol of
$\Sal$, restriction theorem for the Fourier transform and the
Riesz-Thorin interpolation theorem.

\setcounter{equation}{0}
\setcounter{theorem}{0}
\setcounter{lemma}{0}
\section{Preliminaries}\label{sect2}

%\subsection{Notation:}

\subsection{On the analyticity of an integral depending on a parameter}

\begin{lemma} [see e.g. \cite{Rubin2}]\label{l3.1}
Let $f(x,z)$ be  analytic in $z\in{\cal D}\subset\Bbb C$ for
almost all $x\in\R$. If there is a function $F(x)\in L_1(\R)$ such
that $|f(x,z)|\leq F(x)$ for almost all $x\in\R$ and for all
$z\in\cal D$, then the integral $\intR f(x,z)dx$ is an analytic
function on $\cal D$.
\end{lemma}

\subsection{Modified theorem for analytic families of operators}

Let $E$ be the set of simple functions, that is, the set of linear
combinations of characteristic functions of sets in $\R$ having
finite measure. We recall ( see e.g. \cite{Stein2}) , that the
family of operators $\{T_z\}_{z\in S}$, $S=\{z\in\Bbb C: 0\leq\Re
z\leq1\}$, taking the set of simple functions from $L_1(\R)$ into
the space of functions measurable on $\R$, is said to be
admissible growth on $\R$ if for arbitrary  $f, g\in E$, the
function $F(z)=\intR(T_zf)(x)g(x)dx$ is integrable on $\R$ and
has the following properties: \\
i)   $F(z)$ is  analytic  in the interior of $S$;\\
ii)  $F(z)$ is continuous on $S$;\\
iii) $\sup\limits_{\gamma\in\Bbb R^1}\exp\{-a|\Im z|)\}\ln F(z)<\infty$
for some $a<\pi$.\\

We need the following modification of Theorem 1$'$ (\cite{BMO}) which is a special case
of Stein's interpolation theorem \cite{Stein2} (cf. \cite{Nogin}).

\begin{theorem}\label{t3.2}
Let $\{T_z\}_{z\in S}$, $S=\{z\in\Bbb C: 0\leq\Re z\leq1\}$ be an
admissible growth family of multiplier operators satisfying, for
each  $f\in E$, the relations
$$
\|T_{i\gamma} f\|_2\leq M_1(\gamma)\|f\|_2, \quad \quad \quad
\quad \|T_{1+i\gamma} f\|_\infty\leq M_2(\gamma)\|f\|_1,
$$
where $M_j(\gamma)$, $j=1,2$, are independent of $f$, and
$\sup\limits_{\gamma\in\Bbb R^1}\exp\{-a|\gamma|)\} \ln
M_j(\gamma)<\infty$ for some $a<\pi$. Let also the family
$\{T_z\}_{z\in S}$ satisfy  the additional assumption
\begin{eqnarray*}
%|\widehat{T^*_zT_z f}|\leq C_z|\widehat{T_{2\Re z}f}|\quad
|\mu_z(\xi)|^2\leq \mu_{2\Re z}(\xi)  \quad
 \mbox{(pointwise)} \quad \mbox{if} \quad 0<\Re z<\mu \leq 1/2
\end{eqnarray*}
provided $\quad (\widehat{T_z f})(\xi)=\mu_z(\xi)\hat f (\xi)$.
Here $C_z$ is a non-negative function such that $\log C_z\leq K
\exp(k|\Im z|)$ for some $K>0$ and $k<\pi$.
 Then if for some $0<t<1$ and $1/p_t-1/q_t=t$ the operator $T_t$
is  continuous  from $L_{p_t}$ into some topological space $X$ in
which $L_{q_t}$ is continuously embedded, then it is continuous
from $L_{p_t}$ into $L_{q_t}$, provided
\begin{eqnarray*}
\frac{1+t}{2}\leq\frac{1}{p}\leq \frac{1+2t}{2},\quad for \quad
0<t<\mu,
\end{eqnarray*}
and
\begin{eqnarray*}
\frac{1+t}{2}\leq\frac{1}{p}\leq \frac{1+t-2t\mu} {2(1-\mu)},\quad
for \quad \mu\leq t<1.
\end{eqnarray*}

\end{theorem}

\subsection{Uniform asymptotic expansion for the Bessel function}
\begin{lemma} \label{l3.0}
Let $z\in\Omega=\{z\in\C: |z|>\eta,\quad |arg\;z|<\theta\}$, $\eta>0$,
$\theta\in(0;\pi/2)$. Then
\begin{equation}\label{3.0}
J_\nu(z)=(\frac{\pi
z}{2})^{-1/2}\big[e^{-iz}\big(\sum\limits_{m=0}^{M}C^{(\nu)}_{m,-}z^{-m}
+R^{(\nu)}_{M,-}(z)\big)+e^{iz}\big(\sum\limits_{m=0}^{M}C^{(\nu)}_{m,+}z^{-m}
+R^{(\nu)}_{M,+}(z)\big)\big],
\end{equation}
where remainders $R^{(\nu)}_{M,\pm}(z)$ are analytic in $\Omega$ and
\begin{equation}\label{6}
|R^{(\nu)}_{M,\pm}(z)|\leq C|z|^{-M-1}
\end{equation}
 with $C$ independent of $\nu$ and $z$.
\end{lemma}

\subsection{Estimates for the hypergeometric function}
\begin{lemma}
Let  $|z|<1$, $\Re c\geq -1/2$ and there is a constant $M$ such that
$|\Re a|, |\Re b|, |\Re c|\leq M$. Then
\begin{equation}\label{2.5}
|\F (a,b;c;z)| \leq  C
\left\{
\begin{array}{ll}
(\exp(\pi|\Im\; c|)+\frac{|\Gamma(c)|}{|\Gamma(a)\Gamma(b)|}),& \mbox{if} \quad
 a, \; b\ne -m; \quad \Re (c-a-b)>0\\
\exp(\pi|\Im\; c|),& \mbox{if} \quad a=-m \quad or \quad  b=-m,
\end{array}
\right.
\end{equation}
where $m=0,1,2,\dots$ ; $C=C(M)$.
\end{lemma}
\noindent{\bf Proof}. Let us fix $N>0$ and write
$$
\F (a,b;c;z)=(\sum\limits_{m=0}^{N}+\sum\limits_{m=N+1}^{\infty})\frac{(a)_m(b)_m z^m}{(c)_m m!}
\quad.
$$
For $m=0,1,2,\dots, N$, in view of the estimate
\begin{equation}\label{2.6}
|\Gamma(\beta+1/2)|^{-1}\leq C \exp(\pi |\Im\;\beta|),\quad \Re\beta\geq -1,
\end{equation}
with $C$ independent of $\beta$ (see \cite{Bateman1}, p. 48) we
get
\begin{equation}\label{2.7}
|\frac{(a)_m(b)_m}{(c)_m m!}|=|\frac{(a)_m(b)_m \Gamma(c)}{ m!
\Gamma(c+m)}|\leq C(M,N) \exp(\pi|\Im\; c|).
\end{equation}
For $m=N+1,N+2,\dots$, the relation
$$
\lim\limits_{m\to\infty}e^{-a\;\ln\; m}\;\frac{\Gamma(a+m)}{\Gamma(m)}=1
$$
(see \cite{Bateman1}, p. 47) yields
\begin{equation}\label{2.8}
|\frac{(a)_m(b)_m}{(c)_m m!}|\leq C(m) \frac{|\Gamma(c)|}{|\Gamma(a)\Gamma(b)|}
m^{\Re (a+b-c-1)},
\end{equation}
where $\lim\limits_{m\to\infty}C(m)=1$. Making use of (\ref{2.7})
and (\ref{2.8}), we arrive at (\ref{2.5}), provided $a$ and $b$
both differ from $0,-1,-2,\dots\quad$. If $a=-m$ or $b=-m$, then
the hypergeometric series is finite and (\ref{2.5}) follows from
(\ref{2.7}).
$\hfill\blacksquare$

\begin{remark}\label{r3.1}
A simple analysis of the proof shows that it is possible to omit
the restriction $\Re (c-a-b)>0$, if  $|z|\leq 1/2$.
\end{remark}

\subsection{Restriction theorem for the Fourier transform}

Let $\cal S$ be the Schwartz class of rapidly decreasing smooth
functions on $\R$,  and let $\hat f(\sigma)$ denote restriction of
the Fourier transform of $f\in \cal S$ to the unit sphere
$S^{n-1}$ in $\R$. We need the following $(\lp, L_2)$ restriction
property.

\begin{theorem} [\cite{Stein1}]\label{t3.1}
Let $f\in\cal S$. Then
\begin{equation}
\left(\int\limits_{S^{n-1}}|\hat f(\sigma)|^{2}d\sigma\right)^{1/2}\leq C_{p}\normf
\label{3.1}
\end{equation}
if and only if $1\leq p \leq p_0$,  $p_{0}=\frac{2n+2}{n+3}$.
\end{theorem}

\begin{remark}\label{r3.1}
Since $\cal S$ is dense in $\lp$, we can define $\hat f$ on
$S^{n-1}$ for each $f\in\lp$  whenever $1\leq p \leq p_0$.
\end{remark}

\begin{remark}\label{r3.2}
Observe that (\ref{3.1}) holds uniformly  for any sphere of radius
 $1/2\leq r\leq 2$ on place of $S^{n-1}$ .
\end{remark}

\subsection{$(\lp,\lr)$ estimates for certain oscillatory operator}

Let $\psi$ be a fixed smooth function of compact support on $\R$
that vanishes in a neighborhood of the origin, and set

\begin{equation}\label{3.2.1}
(G_{\lambda}f)(x)=\intR \exp(i\lambda|x-y|)\psi(x-y)f(y)dy.
\end{equation}

\begin{lemma}[\cite{Stein1}]\label{l5.1}
Let $1\leq p\leq\frac{2n+2}{n+3}$ for $n>2$, and $1\leq
p<\frac{4}{3}$ for $n=2$. Then
\begin{equation}\label{3.2}
\|G_{\lambda}f\|_{p}\leq C\lambda^{-n/p'}\normf .
\end{equation}
\end{lemma}

\begin{lemma}[ \cite{Bak}, \cite {S}]\label{l5.2}
Let $1\leq p\leq 2$, $q=\frac{(n+1)p'}{n-1}$ for $n>2$, and $1\leq
p<4$, $q=3p'$ for $n=2$. Then
\begin{equation}\label{3.3}
\|G_{\lambda}f\|_{q}\leq C\lambda^{-n/q}\normf.
\end{equation}

\end{lemma}

\setcounter{equation}{0}
\setcounter{theorem}{0}
\setcounter{lemma}{0}
\section{Auxiliary statements}\label{sect3}

\subsection{On the symbol of the operator $\K$}
 We first show that  $\K$ is Fourier multiplier operator.
Set

\begin{equation}
\ma=
\left\{
\begin{array}{clc}
&
\frac{2\pi^{n/2}\Gamma(\alpha)\exp(\alpha \pi i /2)}{\Gamma(n/2)}
\F(\frac{\alpha}{2}, \frac{\alpha+1}{2};\frac{n}{2}; |\xi|^2),
&\mbox{if} \quad |\xi|<1,
\\[14pt]
&
\frac{2^\alpha\pi^{n/2}\Gamma(n/2)}{\Gamma((n-\alpha)/2)}|\xi|^{-\alpha}
\F(\frac{\alpha}{2}, \frac{\alpha-n+2}{2};\frac{1}{2};\frac{1}{|\xi|^2})\>+
\\[14pt] +&
i\frac{2^{\alpha+1}\pi^{n/2}\Gamma((n+1)/2)}{\Gamma((n-\alpha-1)/2)}|\xi|^{-\alpha-1}
\F(\frac{\alpha+1}{2}, \frac{\alpha-n+3}{2};\frac{3}{2};\frac{1}{|\xi|^2}),\quad
&\mbox{if} \quad |\xi|>1.
\end{array}
\right.
\label{4.1}
\end{equation}
%Here, $\F(a,b;c;z)$ is the hypergeometric function.

\begin{lemma}\label{l4.1}
Let $f\in  E$, $0< \Re \alpha<\frac{n}{2}$. Then

\begin{equation}
({\cal F}\K f)(\xi)=\ma({\cal F}f)(\xi)
\label{4.2}
\end{equation}
where the Fourier transform is interpreted as
$$
({\cal F}\K f)(\xi)=\lim\limits_{k\to\infty}^{(L_2)}\int\limits_{|x|<k}(\K f)(x)\exp(ix\cdot\xi)dx.
$$
\end{lemma}
\noindent{\bf Proof}. Let $\chi(|y|)$ be the indicator  of the
unit ball and let $\kf(y)=|y|^{\alpha -n}\exp(i|y|)\chi(|y|)$,
$\ks(y)=|y|^{\alpha -n}\exp(i|y|)(1-\chi(|y|))$. Evidently,
$\kf\in L_{1}(\R)$ for $\Re\alpha>0$, and $\ks\in L_{2}(\R)$ for
$\Re\alpha<n/2$. We put
$$
{\cal K}_{\alpha, N}^{(2)}(y)=
\left\{
\begin{array}{lc}
{\cal K}_{\alpha}^{(2)}(y),& \mbox{if} \quad |y|<N \\
0,& \mbox{if} \quad |y|\geq N.
\end{array}
\right.
$$
Then there is a subsequence ${\cal F K}_{\alpha, N_j}^{(2)}$ such
that $\lim\limits_{j\to\infty}({\cal F K}_{\alpha,
N_j}^{(2)})(\xi)=({\cal F}\ks)(\xi)$ exists almost everywhere.
Therefore (\ref{4.2}) is valid, where
$$
\ma=\lim\limits_{j\to\infty}^{(a.e)}\int\limits_{|y|<N_{j}}(\kf(y)+\ks(y))\exp(iy\cdot\xi)dy=\\
\lim\limits_{j\to\infty}^{(a.e)}\int\limits_{|y|<N_{j}}|y|^{\alpha-n}\exp(i|y|+iy\cdot\xi)dy.
$$
%for almost every $\xi$.
 Making use of the formula
\begin{equation}\label{5}
\intR \exp(ix\cdot y)\varphi(|y|)dy=\frac{(2\pi)^{n/2}}{|x|^{\frac{n-2}{2}}}
\int\limits_{0}^{\infty}\varphi(\rho)\rho^{n/2}J_{n/2-1}(\rho|x|)d\rho,
\end{equation}
(see e.g. \cite{SW}), we obtain
\begin{equation}
\ma=\frac{(2\pi)^{n/2}}{|\xi|^{(n-2)/2}}\int\limits_{o}^{\infty}\rho^{\alpha-n/2}
\exp(i\rho)J_{\frac{n-2}{2}}(\rho|\xi|)d\rho, \label{4.3}
\end{equation}
where the integral is understood as improper  integral provided $
\Re\alpha\geq \frac{n-1}{2}$. Note, that $\ma$ is defined when
$\Re\alpha<\frac{n+1}{2}$.  Now the representation of the integral
in (\ref{4.3}) via hypergeometric functions (see e.g. 2.12.15.3,
\cite{Pr2}) yields (\ref{4.1}). $\hfill\blacksquare$

\begin{remark}\label{r4.1}
%$\lim\limits_{N\to\infty}\int\limits_{o}^{N}\rho^{\alpha-n/2}
%\exp(i\rho)J_{\frac{n-2}{2}}(\rho|\xi|)d\rho$.
In what follows we shall need sometimes to understand (\ref{4.3})
as Abel summable to $\ma$ when $\frac{n-1}{2}\leq
\Re\alpha<\frac{n+1}{2}$, that is,
\begin{equation}\label{4.3.1}
\ma=\lim\limits_{\varepsilon\to 0} \me,
\end{equation}
where
$$
\me=\frac{(2\pi)^{n/2}}{|\xi|^{(n-2)/2}}\int\limits_{o}^{\infty}\rho^{\alpha-n/2}
\exp(i\rho-\varepsilon\rho)J_{\frac{n-2}{2}}(\rho|\xi|)d\rho.
$$
\end{remark}

\begin{lemma}\label{l4.1}
The function $\ma$ admits the following estimates:\\
 1. If ~$0<\Re\alpha\leq\delta_{1}<\frac{n-1}{2}$, then
\begin{equation}\label{4.4}
 \quad |\ma|\leq C \exp(\frac{\pi}{2}|\Im\alpha|),
\end{equation}
where $C=C(n, \delta_{1})$;\\
2. If ~$\frac{n-1}{2}\leq\Re\alpha\leq\delta_{2}<\frac{n+1}{2}$,~
$\alpha\ne \frac{n-1}{2}$, then
\begin{equation}\label{4.5}
\quad |\ma|\leq C \exp(\pi|\Im\alpha|)
\left\{
\begin{array}{lc}
1, & if \quad |\xi|\leq\frac{1}{2} \quad or \quad |\xi|\geq 2, \\
|\Gamma(\alpha-\frac{n-1}{2})|(1+|1-|\xi||^{\frac{n-1}{2}-\Re\alpha}),
& if \quad \frac{1}{2}<|\xi|<2, \quad |\xi|\ne 1
\end{array}
\right.
\end{equation}
 where $C=C(n, \delta_{2})$.\\
3. If ~$\alpha=\frac{n-1}{2}$, then
\begin{equation}\label{4.6}
\quad |\ma|\leq C
\left\{
\begin{array}{lc}
\exp(\frac{\pi}{2}|\Im\alpha|), & if \quad |\xi|\leq\frac{1}{2} \quad or \quad |\xi|\geq 2,\\
1+\ln|1-|\xi||,
& if \quad \frac{1}{2}<|\xi|<2, \quad |\xi|\ne 1,
\end{array}
\right.
\end{equation}
 where $C=C(n)$.
\end{lemma}
\noindent{\bf Proof}. Estimates (\ref{4.4}), provided $|\xi|\ne 1$
and (\ref{4.5}), (\ref{4.6}), provided $|\xi|\leq\frac{1}{2}$ or
$|\xi|\geq 2$ can be readily verified by means of (\ref{2.5}). In
the case $\frac{n-1}{2}\leq\Re\alpha\leq\delta_{2}<\frac{n+1}{2}$,
$\frac{1}{2}<|\xi|<2, \quad |\xi|\ne 1$ we will use the
representation (\ref{4.3.1})
 of the symbol $\ma$. Let $\psi(\rho)$ be a smooth function such that
$\psi(\rho)=0$ if $\rho<1/2$ and $\psi(\rho)=1$ if $\rho>1$. Then
by means of  (\ref{3.0})  we obtain
\begin{equation}\label{4.1.5}
\begin{array}{lll}
\me&=&
\sum\limits_{m=0}^{M}\frac{1}{|\xi|^{\frac{n-1}{2}+m}}\;[C_{m,+}\;
I^{\alpha,\varepsilon}_{m,+}(|\xi|)+C_{m,-}\;
I^{\alpha,\varepsilon}_{m,-}(|\xi|)]\\[14pt]
&+&
\frac{(2\pi)^{n/2}}{|\xi|^{(n-2)/2}}\;I^{\alpha,\varepsilon}(|\xi|)+
\frac{1}{|\xi|^{\frac{n-1}{2}}}\;[C_{M,+}\;
I^{\alpha,\varepsilon}_{M,+}(|\xi|)+C_{M,-}\;
I^{\alpha,\varepsilon}_{M,-}(|\xi|)],
\end{array}
\end{equation}
where
\begin{equation*}
\begin{array}{lll}
I^{\alpha,\varepsilon}(|\xi|)=\int\limits_{0}^{\infty}(1-\psi(\rho))
\rho^{\alpha-n/2}\exp(i\rho-\varepsilon\rho)J_{\frac{n-2}{2}}(\rho|\xi|)d\rho,\\[14pt]
I^{\alpha,\varepsilon}_{m,\pm}(|\xi|)=\int\limits_{0}^{\infty}\psi(\rho)
\rho^{\alpha-\frac{n+1}{2}-m}\exp(i\rho(1\pm |\xi|)-\varepsilon\rho)d\rho,\\[14pt]
I^{\alpha,\varepsilon}_{M,\pm}(|\xi|)=\int\limits_{0}^{\infty}\psi(\rho)
\rho^{\alpha-\frac{n+1}{2}}\exp(i\rho(1\pm |\xi|)-\varepsilon\rho)R_{M,\pm}d\rho.
\end{array}
\end{equation*}
We need only to consider $I^{\alpha,\varepsilon}_{0,\pm}$, since
$|\lim\limits_{\varepsilon\to 0}I^{\alpha,\varepsilon}(|\xi|)|\leq
C(n,\delta_2)$, $|\lim\limits_{\varepsilon\to
0}I^{\alpha,\varepsilon}_{M,\pm}(|\xi|)\leq C(n,\delta_2)$, and
the integrals $I^{\alpha,\varepsilon}_{m,\pm}$, $m=1,2,\dots$,
have better decay properties then
$I^{\alpha,\varepsilon}_{0,\pm}$. The calculation in the case
$\frac{n-1}{2}<\Re\alpha\leq\delta_{2}$ is easy by virtue of
(2.3.3.1), \cite{Pr1}:
$$
I^{\alpha,\varepsilon}_{0,\pm}(|\xi|)=\int\limits_{0}^{\infty}(\psi(\rho)-1)
\rho^{\alpha-\frac{n+1}{2}}\exp(i\rho(1\pm|\xi|)-\varepsilon\rho)d\rho+
\Gamma(\alpha-\frac{n-1}{2})[i(1\pm|\xi|)-\varepsilon]^{\frac{n-1}{2}-\alpha}.
$$
Thus,
$$
|\lim\limits_{\varepsilon\to 0}\me|\leq C \exp(\pi|\Im\alpha|)
|\Gamma(\alpha-\frac{n-1}{2})|(1+|1-|\xi||^{\frac{n-1}{2}-\Re\alpha}).
$$
Then, letting $\alpha=\frac{n-1}{2}+i\gamma$, $\gamma\in\Bbb R^1$
and  integrating by parts yields
\begin{equation*}
\begin{array}{lll}
I^{\alpha,\varepsilon}_{0,\pm}(|\xi|)&=&
-\frac{i(1+|\xi|)-\varepsilon}{i\gamma}
\int\limits_{0}^{\infty}(\psi(\rho)-1)
\rho^{i\gamma}\exp(i\rho(1\pm|\xi|)-\varepsilon\rho)d\rho\\[14pt]
&-&
\Gamma(i\gamma)[i(1\pm|\xi|)-\varepsilon]^{-i\gamma}-
\int\limits_{0}^{\infty}\psi'(\rho)
\rho^{i\gamma}\exp(i\rho(1\pm|\xi|)-\varepsilon\rho)d\rho,
\end{array}
\end{equation*}
if $\gamma\ne 0$, and
\begin{equation*}
\begin{array}{lll}
I^{\alpha,\varepsilon}_{0,\pm}(|\xi|)&=&
c+\ln(\varepsilon-i(1\pm|\xi|))-
\int\limits_{0}^{\infty}(\psi(\rho)-1)
\ln\rho \exp(i\rho(1\pm|\xi|)-\varepsilon\rho)d\rho-\\[14pt]
&-&
\int\limits_{0}^{\infty}\psi'(\rho)
\ln\rho \exp(i\rho(1\pm|\xi|)-\varepsilon\rho)d\rho,
\end{array}
\end{equation*}
if $\gamma=0$, in view of 4.331.1 \cite{Gr}, where $c$ is Euler's
constant.  The desired result then readily follows from this.
$\hfill\blacksquare$

\subsection{Auxiliary estimates}

Let us consider  the operator $\Sa$ defined by (\ref{4.6.1}).
Decompose $\Sa$ into (\ref{4.7}) and denote by ${\cal
K}^\alpha_\ell$ the kernel of the operator $S^\alpha_\ell$.

\begin{lemma}\label{l4.2}
Let $\ell=0,1,2,\dots$. Then
\begin{equation}\label{4.8}
 \quad |\hat{{\cal K}^\alpha_\ell}(\xi)|\leq C2^{-M\ell} \left\{
\begin{array}{ll}
1, & if \quad |\xi|\leq\frac{1}{2}   , \\
(1+|\xi|)^{-M},
& if \quad |\xi|\geq 2
\end{array}
\right.
\end{equation}
for any $M>0$, and
\begin{equation}\label{4.9}
 \quad |\hat{{\cal K}^\alpha_\ell}(\xi)|\leq
C2^{(\alpha-\frac{n-1}{2}) \ell}, \quad if \quad
\frac{1}{2}<|\xi|<2.
\end{equation}

\end{lemma}
\noindent{\bf Proof}. Let $\ell=1,2,\dots$. Making use of
(\ref{5}), we get
$$
\hat{{\cal K}^\alpha_\ell}(\xi)=
\frac{2^{\ell(\alpha-n/2+1)}(2\pi)^{n/2}}{|\xi|^{\frac{n-2}{2}}}
\int\limits_{0}^{\infty}\rho^{n/2}\exp(i2^\ell\rho)J_{\frac{n-2}{2}}(2^\ell\rho\xi)\psi(\rho)d\rho.
$$
In the case $|\xi|\leq\frac{1}{2}$ we apply the integral representation for the Bessel
function
$$
J_\nu(z)=\frac{(z/2)^\nu}{\sqrt\pi\Gamma(\nu+\frac{1}{2})}
\int\limits_{-1}^{1}(1-t^2)^{\nu-1/2}e^{izt}dt,
$$
(see \cite{Bateman2}) and obtain
$$
\hat{{\cal K}^\alpha_\ell}(\xi)=C(n)2^{\alpha\ell}
\int\limits_{-1}^{1}(1-t^2)^{\frac{n-3}{2}}dt
\int\limits_{0}^{\infty}\rho^{n/2}\exp(i2^\ell\rho+i2^\ell\rho\xi t)\psi(\rho)d\rho.
$$
Then the estimate (\ref{4.8}) for $|\xi|\leq\frac{1}{2}$ can be easily obtained
 from the relation
$$
\int\limits_{0}^{\infty}\rho^{n/2}\exp(i2^\ell\rho+i2^\ell\rho\xi
t)\psi(\rho)d\rho= \Big[\frac{i}{2^\ell(1+\xi
t)}\Big]^k\int\limits_{0}^{\infty} \big(\frac{d}{d\rho}\big)^k
(\rho^{n/2}\psi(\rho))\exp(i2^\ell\rho+i2^\ell\rho\xi t)d\rho.
$$
The latter obviously holds for any $k>0$ by means of repeated integration by parts.\\
% that
%$|\hat{{\cal K}^\alpha_\ell}(\xi)|\leq C2^{\ell(\alpha-k)}$ for any $k>0$ and
%$|\xi|\leq\frac{1}{2}$ with $C$ independent on $\ell$.
Letting $|\xi|>\frac{1}{2}$, in view of  (\ref{3.0})  we have
\begin{eqnarray*}
\hat{{\cal K}^\alpha_\ell}(\xi)= 2^{\ell(\alpha-\frac{n-1}{2})}
\big[\sum\limits_{m=0}^{M-1}\frac{1}{|\xi|^{\frac{n-1}{2}-m}}
(C_{m,-}I_{m,-}(|\xi|)+C_{m,+}I_{m,+}(|\xi|)) + \\
\nonumber \frac{1}{|\xi|^{\frac{n-1}{2}}}(C_{M,-}I_{M,-}(|\xi|)+
C_{M,+}I_{M,+}(|\xi|)\big],
\end{eqnarray*}
where
\begin{equation*}
\begin{array}{ll}
I_{m,\pm}(|\xi|)=
\int\limits_{0}^{\infty}\rho^{\frac{n-1}{2}-m}\exp(i2^\ell\rho(1\pm|\xi|))\psi(\rho)d\rho,\\[14pt]
I_{M,\pm}(|\xi|)=\int\limits_{0}^{\infty}\rho^{\frac{n-1}{2}}\exp(i2^\ell\rho(1\pm|\xi|))\psi(\rho)
R_{M,\pm}(2^\ell\rho |\xi|)d\rho.
\end{array}
\end{equation*}
Evidently, $|\hat{{\cal K}^\alpha_\ell}(\xi)|\leq C2^{\ell(\alpha-\frac{n-1}{2})}$, provided
$\frac{1}{2}<|\xi|< 2$.
Repeated integration by parts yields
$$
\int\limits_{0}^{\infty}\rho^{\frac{n-1}{2}-m}\exp(i2^\ell\rho(1\pm|\xi|))\psi(\rho)d\rho=
\Big[\frac{i}{2^\ell(1\pm\xi)}\Big]^k\int\limits_{0}^{\infty}
\big(\frac{d}{d\rho}\big)^k
(\rho^{\frac{n-1}{2}-m}\psi(\rho))\exp(i2^\ell\rho(1\pm |\xi|
))d\rho.
$$
 Taking into account (\ref{6}), we obtain (\ref{4.8}) for
$|\xi|\geq 2$.
 $\hfill\blacksquare$

%$$
%|\hat{{\cal K}^\alpha_\ell}(\xi)|\leq C2^{-M\ell}(1+|\xi|)^{-M},\quad |\xi|\geq 2.
%$$

\begin{lemma}\label{l4.3}
Let $\ell=0,1,2,\dots$, $f\in\lp$, $1\leq p\leq\frac{2(n+1)}{n+3}$. Then
\begin{equation}\label{4.9}
\|\Sa_{\ell}f\|_2\leq C2^{\ell(\alpha-n/2)}\normf
\end{equation}
with C independent of $\ell$.
\end{lemma}
\noindent{\bf Proof}. Since the kernel ${\cal K}^\alpha_\ell$ of the operator $\Sal$ is supported
in the annulus $2^{\ell-1}<|y|<2^{\ell+1}$, it suffices to proof (\ref{4.9}) for functions supported
in the ball of radius $2^{\ell+1}$. For such $f$ we obtain
$$
\|\Sa_{\ell}f\|^2_2=\Big(\int\limits_{|\xi|\leq
1/2}+\int\limits_{|\xi|\geq 2}+
\int\limits_{1/2<|\xi|<2}\Big)|\hat{{\cal K}^\alpha_\ell}(\xi)|^2  |\hat f(\xi)|^2 d\xi=\\
I_1+I_2+I_3.
$$
The integrals $I_1$ and  $I_2$ are easily treated by  (\ref{4.8}).
% $\hat{\cal K}^\alpha_\ell$ for $|\xi|\leq\frac{1}{2}$ and $|\xi|\geq 2$.
 For  $I_3$ we use  (\ref{3.1}), the radiality
of
 $\hat{\cal K}^\alpha_\ell$ and the estimate
 $\|{\cal K}^\alpha_\ell\|_2\leq C2^{\ell(\alpha-n/2)}$.
Hence
$$
I_3\leq \normf^2 \|{\cal K}^\alpha_\ell\|_2^2\leq
C2^{2\ell(\alpha-n/2)}\normf^2,
$$
and (\ref{4.9}) is proved.
$\hfill\blacksquare$

\subsection{$\lp-$ boundedness property of the operator $\Sa$}\label{sect4}

\begin{theorem}\label{t2}
The operator $\Sa$ is bounded on $\lp$ whenever
\\
1.  $0<\alpha<\frac{1}{2}$, $\frac{2}{2-\alpha}<p<\frac{2}{\alpha}$, if \quad $n=2$;\\
2.  $0<\alpha<\frac{n(n-1)}{2(n+1)}$, $\frac{n}{n-\alpha}<p<\frac{n}{\alpha}$, if \quad $n\geq 3$;\\
3.  $\frac{n(n-1)}{2(n+1)}\leq\alpha<\frac{n-1}{2}$,
$\frac{2(n-1)}{3(n-1)-4\alpha}<p<\frac{2(n-1)}{4\alpha-n+1}$, if \quad $n\geq 3$.
\end{theorem}

\noindent{\bf Proof}.
%We recall (\ref{0.2}). Since $\N$ is bounded
%on $\lp$ for any $p$, it remains to consider  $\Sa$.
Making use of
(\ref{new:4.2}), we have
 $\|\Sa_\ell\|=2^{\alpha\ell}\|G_{2^{\ell}}\|$, and
\begin{equation}
\|\Sa_{\ell}f\|_{p}\leq C2^{\ell(\alpha-\frac{n}{p'})}\normf,
\label{5.1}
\end{equation}
for $1\leq p\leq\frac{2n+2}{n+3}$  when $n>2$ and $1\leq
p<\frac{4}{3}$ when $n=2$ in view of Lemma~\ref{l5.1}. Summation
over $\ell$ and duality give the desired conclusions in the first
two cases.\\
To obtain the last desired result, we will interpolate between
(\ref{5.1}) and the estimate
\begin{equation}\label{4.1.1}
\|\Sa_{\ell}f\|_{2}\leq C2^{\ell(\alpha-\frac{n-1}{2})}\|f\|_2.
\end{equation}
By this we arrive at the inequality
$$
\|\Sa_{\ell}f\|_{p}\leq C2^{\ell(\alpha+\frac{n-1}{2p}-\frac{3(n-1)}{4})}\normf,
$$
and we are done.
$\hfill\blacksquare$

\begin{remark}
A necessary condition for $\Sa$ to be bounded on $\lp$ is
$\frac{n}{n-\alpha}<p<\frac{n}{\alpha}$.
\end{remark}

\subsection{$(\lp, \lr)$ estimates along a segment through $(\frac{1}{p}, \frac{1}{p'})$
and perpendicular to the line of duality}\label{s4.4}

\begin{lemma}\label{l4.4.1}
Let $n=2$, $\frac{1}{2}<\alpha <2$, $f\in\lp$ and   $(\frac{1}{p}, \frac{1}{q})$
is  on the open segment  $(B'B)$.
 Then there is $C=C(\alpha)$ such that $\|\Sa f\|_q\leq C\normf$.
\end{lemma}
\noindent{\bf Proof}.
Let us fix point $P(\frac{1}{p};\frac{1}{q})$ on $(BB')$ such that
$\frac{1}{3}-\frac{\alpha}{6}<\frac{1}{q}<\min\{\frac{1}{4},1-\frac{\alpha}{2}\}$.
Making use of (\ref{3.3}) and the assertion
\begin{equation}\label{new:5.1}
(\Sa_{\ell}f)(x)=2^{\alpha\ell}G_{2^\ell}(f_{2^\ell})(2^{-\ell}x)
\end{equation}
where $f_{\lambda}(x)=f(\lambda x)$, we have
\begin{equation}\label{5.2}
\|\Sa_{\ell}f\|_q\leq C2^{(\alpha-2+6/q)\ell}\|f\|_r,\quad q=3r'.
\end{equation}
Note that by duality (\ref{3.3}) is equivalent to the estimate
\begin{equation}\label{5.3}
\|G_{\lambda}f\|_{b}\leq C\lambda^{-2/3b}\|f\|_a,
\end{equation}
for  $b>4/3$ and $1/b=3(1-1/a)$. It follows from (\ref{new:5.1}) and (\ref{5.3}) that
\begin{equation}\label{5.4}
\|\Sa_{\ell}f\|_q\leq C2^{(\alpha-2+2/q)\ell}\|f\|_a,
\end{equation}
where $1/a=1-1/3q$. Estimates (\ref{5.2}) and (\ref{5.4}) applied to
a characteristic function $f=\chi_E$ give
\begin{equation*}
%\label{5.5}
\|\Sa_{\ell}\chi_E\|_q\leq C\min \{2^{(\alpha-2+6/q)\ell}|E|^{1-3/q},\quad
2^{(\alpha-2+2/q)\ell}|E|^{1-1/3q}\}.
\end{equation*}
Observe that the first term in the braces is smaller then the second
precisely when $2^\ell<|E|^{2/3}$.
Let $v=|E|$ and for $v>0$ let $N=N(v)$ be the integer such that
$2^N<v^{2/3}\leq 2^{N+1}$. Then

\begin{equation*}
\sum\limits_{\ell=1}^{\infty}\|\Sa_{\ell}\chi_E\|_q
\leq C \sum\limits_{\ell=-\infty}^{N}2^{(\alpha-2+6/q)\ell}v^{1-3/q}+
C \sum\limits_{\ell=N+1}^{\infty}2^{(\alpha-2+2/q)\ell}v^{1-1/3q}.
\end{equation*}
Since $\frac{1}{3}-\frac{\alpha}{6}<\frac{1}{q}<1-\frac{\alpha}{2}$,
the two geometric series are convergent, and we obtain
 \begin{eqnarray*}
\sum\limits_{\ell=1}^{\infty}\|\Sa_{\ell}\chi_E\|_q
\leq C 2^{N(\alpha-2+6/q)\ell}v^{1-3/q}+
C 2^{(N+1)(\alpha-2+2/q)\ell}v^{1-1/3q}\leq\\ \nonumber
\leq C v^{(2/3)(\alpha-2+6/q)\ell}v^{1-3/q}+v^{(2/3)(\alpha-2+2/q)\ell}v^{1-1/3q}
\leq C v^{1/q-1/3+2\alpha/3}.
\end{eqnarray*}
 Hence we get

\begin{equation}\label{5.7}
\|\Sa_{\ell}\chi_E\|_q\leq C\|\chi_E\|_p,
\end{equation}
whenever
$\frac{1}{3}-\frac{\alpha}{6}<\frac{1}{q}<\min\{\frac{1}{4};1-\frac{\alpha}{2}\}$
and $\frac{1}{p}=\frac{1}{q}-\frac{1}{3}+\frac{2\alpha}{3}$. It is
well known (\cite{SW}) that (\ref{5.7}) is equivalent to the
Lorentz space estimate
\begin{equation}\label{5.8}
\|\Sa f\|_q\leq C\|f\|_{L^{p,1}}
\end{equation}
for the same $(p,q)$. In the case $\alpha\geq 3/2$ the estimate (\ref{5.8})
implies
\begin{equation}\label{5.8.1}
\|\Sa f\|_q\leq C\normf
\end{equation}
 for any given point $(1/p,1/q)$ on $(BB')$. For $\alpha<3/2$
the estimate (\ref{5.8.1}) for any given point
$(\frac{1}{p},\frac{1}{q})$ on $(BB')$ follows by duality and
interpolation.
$\hfill\blacksquare$

\begin{lemma}\label{l4.4}
Let $n>2$, $\frac{n-1}{2}<\alpha <n$, $f\in\lp$ and   $(\frac{1}{p}, \frac{1}{q})$
is either on the open segment $(P'P)$, provided $\frac{n-1}{2}<\alpha<\frac{n}{2}$ or on $(B'B)$,
provided $\frac{n}{2}\leq\alpha<n$.
 Then there is $C=C(\alpha)$ such that $\|\K f\|_q\leq C\normf$.
\end{lemma}
\noindent{\bf Proof}. Let us consider the family of operators
$T_z$, $z\in S=\{z\in \Bbb C : 0\leq Rez \leq 1\}$, defined on simple
functions $f$ as follows:
\begin{eqnarray}
T_z f=&
\left\{
\begin{array}{lc}
\frac{1}{\Gamma(\frac{z(n+1)}{2})}K^{\frac{z(n+1)}{2}+\frac{n-1}{2}}f,& \mbox{if} \quad z\ne 0\\
0,& \mbox{if}\quad z=0.
\end{array}
\right.
\label{new:4.1}
\end{eqnarray}
We prove that this family is admissible growth. Given $f, g \in
E$, let us put $F(z)=\intR (T_z f)(x)g(x)dx$. The condition i)
obviously holds by virtue of Lemma~\ref{l3.1}. To prove ii) it
remains to show that $F(z)$ is continuous in a strip
$S_\delta=\{z\in\C:0\leq \Re z\leq\delta\}$ of arbitrarily small
width $\delta$. Let us fix $\delta<\frac{1}{n+1}$. Since $T_z f\in
L_2$ for $z\in S_\delta$ the application of Parseval's formula
yields
\begin{equation}
F(z)=\frac{(2\pi)^{-n}}{\Gamma(\frac{z(n+1)}{2})}\intR m_{\frac{z(n+1)}{2}+\frac{n-1}{2}}
(|\xi|)(Ff)(\xi)(Fg)(-\xi)d\xi.
\label{4.10}
\end{equation}
The integral on the right-hand side is continuous on $S_\delta$, provided $z\ne 0$,
in view of (\ref{4.5}). In the case $z=0$ we have
$\lim\limits_{z\to 0} T_z f =0$, hence $F(z)$ is continuous at the point $z=0$.
Assertion iii) easily follows from (\ref{4.10}) and (\ref{4.5}), if
 $0\leq \Re z \leq\delta$, $\delta<\frac{1}{n+1}$.
For $\frac{1}{n+1}<Re\; z\leq 1$ we have

\begin{equation}\label{4.11}
|F(z)|\leq\intR|g(x)|(|Q_z(x)|+|R_z(x)|)dx,
\end{equation}
where
$$
Q_z(x)=\frac{1}{\Gamma(\frac{z(n+1)}{2})}
\int\limits_{|y|< 1}\frac{\exp(i|y|)}{|y|^{\frac{n+1}{2}(1-z)}}f(x-y)dy
$$
and
$$
R_z(x)=\frac{1}{\Gamma(\frac{z(n+1)}{2})}
\int\limits_{|y|>1}\frac{\exp(i|y|)}{|y|^{\frac{n+1}{2}(1-z)}}f(x-y)dy.
$$
Then iii) follows from the estimates
$$
|Q_z(x)|\leq \frac{1}{|\Gamma(\frac{z(n+1)}{2})|}
\int\limits_{|y|< 1}\frac{|f(x-y)|}{|y|^{\frac{n}{2}}}dy\leq
A_{n}(x) \exp\big (\frac{\pi (n+1)}{2}|\Im\; z|\big )
$$
and
$$
|R_z(x)|\leq\frac{1}{|\Gamma(\frac{z(n+1)}{2})|}
\int\limits_{|y|>1}|f(x-y)|dy \leq C_n \exp\big (\frac{\pi
(n+1)}{2}|\Im\; z|\big ),
$$
where $A_n(x)\in\lp$, $1\leq p\leq\infty$. We apply
Theorem~\ref{t3.2} to the family (\ref{new:4.1}). For each
$\gamma\in{\Bbb R}^1$ it follows from (\ref{4.5}) and (\ref{4.10})
that
$$
\|T_{i\gamma}f\|_2\leq C \exp\big(\frac{3\pi
(n+1)}{4}|\gamma|\big) \|f\|_2,
$$
$$
\|T_{1+i\gamma}f\|_\infty\leq C \exp\big(\frac{\pi
(n+1)}{2}|\gamma|\big) \|f\|_1.
$$
Moreover,
\begin{equation}\label{new:4.10}
|\mu_z(|\xi|)|^2\leq C \frac{\exp(\frac{\pi\Im
z(n+1)}{2})\Gamma(\Re z(n+1))} {|\Gamma(\frac{z(n+1)}{2})|^2}\;
\mu_{2\Re z}(|\xi|),\quad 0<\Re z<\frac{1}{n+1},
\end{equation}
where
$\mu_z(|\xi|)=\frac{1}{\Gamma(\frac{z(n+1)}{2})}m_{\frac{z(n+1)}{2}+\frac{n-1}{2}}(|\xi|)$
(see Appendix). Next, we show that the operator $T_t$, $t\in
[0;1]$, is continuous from $L_{p_t}$   into $X=L_{1,\beta}= \{
f:\intR\frac{|f(x)|}{(1+|x|)^\beta}dx<\infty,\quad \beta>n\}$,
provided
 $\frac{1}{p_t}>1-\frac{(n+1)(1-t)}{2n}$.
Let $f\in L_{p_t}$,
then
$$
|(T_tf)(x)|\leq \frac{1}{|\Gamma(\frac{z(n+1)}{2})|}(
\int\limits_{|y|< 1}+\int\limits_{|y|>1})\frac{|f(x-y)|}{|y|^{\frac{n+1}{2}(1-t)}}dy=
\frac{1}{|\Gamma(\frac{z(n+1)}{2})|}\{(T_{o,t}f)(x)+(T_{\infty,t}f)(x)\}.
$$
The operator $T_{o,t}$ is bounded on  $L_{p_t}$, while the
operator $T_{\infty,t}$ is bounded from $L_{p_t}$ into $L_\infty$.
Putting $t=\frac{2\alpha-n+1}{n+1}$,
$\frac{1}{p_t}-\frac{1}{q_t}=t$ and interpolating, we obtain the
desired result.
$\hfill\blacksquare$

\setcounter{equation}{0}
\setcounter{theorem}{0}
\setcounter{lemma}{0}
\section{Proof of the main result}\label{sect4}

We now turn to the proof of Theorem~\ref{t1}

I. Let us decompose $\K$ into (\ref{0.2}). Since the operator
$N^\alpha$ is  bounded  from $\lp$ into $\lr$, provided
$0<\alpha<n$, $1<p<\frac{n}{\alpha}$,
 $\frac{1}{q}=\frac{1}{p}-\frac{n}{\alpha}$
by Hardy-Littlwood-Sobolev theorem and the mapping $f\to N^\alpha f$ is of
"weak-type" $(1, 1-\frac{\alpha}{n})$,
we obtain
$[OFAA']\setminus(\{A\}\cup \{A'\})\subset\cal L(N^\alpha)$.
So it remains to investigate the operator $\Sa$. We will consider  the following
situations separately.

1). The cases  $0<\alpha<\frac{n(n-1)}{2(n+1)}$, $n>2$ and  $0<\alpha<\frac{1}{2}$,
$n=2$.

In view of  H\"older's inequality we have  $\|\Sa f\|_\infty\leq C\normf$, provided
$\frac{1}{p}>\frac{\alpha}{n}$. It follows from theorem~\ref{t2} and the Riesz-Thorin
interpolation theorem that $(A'H'HAE)\cup (H'H)\cup (A'E]\cup [EA)\in {\cal L}(\Sa)$.

2). In the cases $1/2<\alpha<2$, $n=2$ and $\frac{n}{2}\leq\alpha<n$, $n>2$
the conclusion $(A'B'BAE)\cup (B'B)\cup (A'E]\cup [EA)\in {\cal L}(\Sa)$
is an easy consequence of Lemma~\ref{l4.4.1} and Lemma~\ref{l4.4}.

3). The case $\frac{n(n-1)}{2(n+1)}<\alpha<\frac{n}{2}$, $\alpha\neq\frac{n-1}{2}$, $n>2$.

Having decomposition (\ref{new:4.2}) of the operator $\Sa$ we shall now apply
Riesz-Thorin interpolation theorem for the operators (\ref{4.7}).
Let us first interpolate between the estimates (\ref{5.1})
 and (\ref{4.9}), taking both of them at the endpoint  $p_0=\frac{2n+2}{n+3}$. Then
we arrive at the inequality
\begin{equation}\label{3.4}
\|\Sal f\|_q\leq C2^{\ell(\frac{n}{q}+\alpha-n)}\|f\|_{p_0},
\end{equation}
which is valid for any $p_0\leq q\leq 2$. We note, that in view of (\ref{3.4}) the operator
$\Sa$ is bounded from $\L_{p_0}$ into $\lr$ in the case $\alpha\leq n/2$ whenever
$\frac{1}{q}<1-\frac{\alpha}{n}$. Our second step will be interpolation between
(\ref{3.4}) and the estimate
\begin{equation}\label{5.12}
\|\Sa_{\ell}f\|_\infty\leq C2^{\ell(\alpha-n)}\|f\|_1,
\end{equation}
which gives
\begin{equation}\label{5.13}
\|\Sal f\|_q\leq C2^{\ell(\frac{n}{q}+\alpha-n)}\normf,
\end{equation}
where  $1\leq p\leq p_0$, $q=\frac{(n-1)p'q_1}{2(n+1)}$, $p_0\leq q_1\leq 2$.
The exponent of $2^\ell$ is negative if
$\frac{1}{q}<1-\frac{\alpha}{n}$. By virtue of the restriction
$\frac{1}{q}\leq\frac{1}{q_1}$ we have to require
 $\frac{1}{q_1}\geq 1-\frac{\alpha}{n}$.
Then $\frac{1}{p}>\frac{n+3}{2(n+1)}$ and (\ref{5.8.1}) is valid
for such $(p,q)$. It follows immediately from Lemma~\ref{l4.4} and
Theorem~\ref{t2}, that
\begin{equation*}
\begin{array}{ll}
(A'G'P'PGAE)\cup(A'E]\cup[EA)\cup (P'P)\subset{\cal L}(\Sa),\quad \mbox{if}
\quad  \frac{n-1}{2}<\alpha<\frac{n}{2},\\
(A'G'CC'GAE)\cup(A'E]\cup[EA)\cup(CC')\subset{\cal L}(\Sa),\quad  \mbox{if}
\quad \frac{n(n-1)}{2(n+1)}\leq\alpha<\frac{n-1}{2}.
\end{array}
\end{equation*}

4). The case $n=2$, $\alpha=\frac{1}{2}$.

Interpolation between (\ref{5.1}) and (\ref{5.12}) yields
$$
 \|\Sal f\|_{q}\leq 2^{\ell(\frac{2}{q}+\alpha-2)}\normf,
$$
where $\frac{1}{p}=1-\frac{1}{q}(p_1-1)$, provided $1\leq p_1<\frac{4}{3}$.
Thus, (\ref{5.8.1}) holds whenever $\frac{1}{p}>1-\frac{2-\alpha}{6}$,
$\frac{1}{q}<1-\frac{\alpha}{2}$. In particular, this is the case for
$\alpha=\frac{1}{2}$ and consequently
$(A'H'H'AE)\cup(A'E]\cup[EA)\subset{\cal L}(\Sa)$.

5). The case $n>2$, $\alpha=\frac{n-1}{2}$.
In a similar way,
making use of the interpolation between (\ref{4.1.1}) and (\ref{5.12}),
we obtain
$$
\|\Sal f\|_{p'}\leq 2^{\ell(\frac{n+1}{p'}+\alpha-n)}\normf.
$$
It is clear that $(DE]\subset{\cal L}(\Sa)$ and therefore
$(A'G'LGAE)\cup(A'E]\cup[EA)\subset{\cal L}(\Sa)$
for $\alpha=\frac{n-1}{2}$ by means of (\ref{5.13}).

II. One can show by an elementary argument that the operator $\K$
is unbounded from $\lp$ into $\lr$ whenever either $(\frac{1}{p},
\frac{1}{q})\in [HAF]\cup[H'A'O]$ or $(\frac{1}{p},
\frac{1}{q})\in [A'AE]\setminus (AA')$. The proof of the negative
result 3) is precisely that used in \cite{Borjeson}.
$\hfill\blacksquare$

\setcounter{equation}{0}
\setcounter{theorem}{0}
\setcounter{lemma}{0}
\section{Appendix}

To verify (\ref{new:4.10}),  we first obtain the auxiliary estimate
\begin{equation}\label{4.1.2}
|1-|\xi|^2|^{-\Re z(n+1)}\leq C_z \mu_{2\Re z}(|\xi|),\quad \Re
z<\frac{1}{n+1},
\end{equation}
where  $C_z=C\Gamma(\Re z(n+1))$ (with $C$ independent of $z$). In
the case $|z|<1$ we have
\begin{eqnarray*}
\mu_{2\Re z}(|\xi|)= \frac{2\pi^{n/2}\Gamma(\Re
z(n+1)+\frac{n-1}{2})\exp[(\Re z(n+1)+\frac{n-1}{2})\frac{\pi
i}{2}]}
{\Gamma(\frac{n}{2})}\quad(1-|\xi|^2)^{-\Re z(n+1)}\\
\times\F\big(\frac{n+1}{4}-\frac{\Re z(n+1)}{2},
\frac{n-1}{4}-\frac{\Re z(n+1)}{2};\frac{n}{2};|\xi|^2\big).
\end{eqnarray*}
Making use of the integral representation for the hypergeometric function
$$
\F(a,b;c;z)=\frac{\Gamma(c)}{\Gamma(b)\Gamma(c-b)}\int\limits_0^1
t^{b-1}(1-t)^{c-b-1}(1-tz)^{-a}dt,  \quad 0<\Re b<\Re c
$$
(see \cite{Pr3}, p. 431), it is not difficult to verify that
$$
\mu_{2\Re z}(|\xi|)\geq C \Gamma(\Re z(n+1)+\frac{n-1}{2})
(1-|\xi|^2)^{-\Re z(n+1)},
$$
and therefore (\ref{4.1.2}) is proved provided $|\xi|<1$. In the
case $|\xi|>1$ we use the representation (\ref{4.3.1}) of the
symbol $m_{\Re z(n+1)+\frac{n-1}{2}}(|\xi|)$. Applying
decomposition (\ref{4.1.5}) and the arguments
 used in the proof
of Lemma~\ref{l4.1}, we obtain (\ref{4.1.2}). Since
\begin{eqnarray*}
|\mu_z(|\xi|)|^2=\frac{4\pi^n|\Gamma(\frac{z(n+1)}{2}+\frac{n-1)}{2})|^2
|\exp[(z(n+1)+\frac{n-1}{2})\frac{\pi i}{2}]|^2}
{|\Gamma(\frac{z(n+1)}{2}|^2 |\Gamma(\frac{n}{2})|^2}(1-|\xi|^2)^{-\Re z(n+1)}\\
\nonumber \times|\F\big(\frac{n+1}{4}-\frac{z(n+1)}{4},
\frac{n-1}{4}-\frac{z(n+1)}{4};\frac{n}{2};|\xi|^2\big)|^2  ,\quad
|\xi|<1,
\end{eqnarray*}
it follows from (\ref{4.1.2}) and (\ref{2.5}) that
(\ref{new:4.10}) is valid. We prove (\ref{new:4.10}) for $|\xi|>1$
in just the  same way as is proved for $|\xi|<1$ .\\[14pt]

 {\bf Acknowledgments.} I would like to express my  gratitude to
Prof. V.\,A.~ Nogin for suggesting this problem to me as well as
for  his patient help and criticisms. I am  thankful to Prof.
E.~Liflyand for numerous discussions and  valuable comments to the
manuscript of the paper, and to Prof. B.~Rubin  for  useful
discussions of the results.

\end{document}